\documentclass[reqno,final]{amsart}
\usepackage{natbib}  
\usepackage{fancyhdr} 
\usepackage{color} 
\usepackage{hyperref} 
\usepackage{graphicx} 

\usepackage{amsthm}
\usepackage{amssymb}
\usepackage{pstricks}
\usepackage{amsmath, amsfonts}
\usepackage{dsfont}

\definecolor{aleacolor}{rgb}{0.16,0.59,0.78}

\hypersetup{
breaklinks,
colorlinks=true,
linkcolor=aleacolor,
urlcolor=aleacolor,
citecolor=aleacolor}

\renewcommand{\headheight}{11pt}
\renewcommand{\cite}{\citet}

\theoremstyle{plain}
\newtheorem{theorem}{Theorem}[section]
\newtheorem{proposition}[theorem]{Proposition}
\newtheorem{lemma}[theorem]{Lemma}
\newtheorem{corollary}[theorem]{Corollary}

\theoremstyle{definition}
\newtheorem{definition}[theorem]{Definition}
\theoremstyle{remark}
\newtheorem{remark}[theorem]{Remark}

\makeatletter \@addtoreset{equation}{section} \makeatother
\renewcommand\theequation{\thesection.\arabic{equation}}
\renewcommand\thefigure{\thesection.\arabic{figure}}
\renewcommand\thetable{\thesection.\arabic{table}}

\newcommand{\aleaIndex}[1]{\href{http://alea.impa.br/english/index_v#1.htm}{\bf #1}}
\eheader{ALEA, Lat. Am. J. Probab. Math. Stat.}{\aleaIndex{14}}{2017}{201}{218}

\newcommand{\N}{\mathbb{Z}_{+}}

\newcommand{\R}{\mathbb{R}}

\begin{document}

\title[RBSDEs when the obstacle is not RC in a general filtration]{Reflected BSDEs when the obstacle is not right-continuous in a general filtration}

\author{Baadi Brahim}
\author{Ouknine Youssef}

\address{Ibn Tofa\"{\i}l University,\newline
Department of mathematics, faculty of sciences,\newline
BP 133, K\'{e}nitra, Morocco }
\email{baadibrahim@gmail.com}

\address{Cadi Ayyad University,\newline
Av. Abdelkrim Khattabi \newline
40000, Gu\'{e}liz Marrakesh, Morocco, and\newline
Hassan II Academy of Sciences and Technology, Rabat, Morocco.}
\email{ouknine@uca.ma}


\subjclass[2000]{60K35, 82B43.}
\keywords{ Backward stochastic differential equation, Reflected backward stochastic differential
equation, General filtration, Strong optional supermartingale,
Mertens decomposition.}

\begin{abstract}
  We prove existence and uniqueness  of the  reflected backward
stochastic differential equation's (RBSDE) solution  with a lower
obstacle which is assumed to be right upper-semicontinuous but not
necessarily right-continuous in a  filtration that supports a
Brownian motion $W$ and an independent Poisson random measure $\pi$.
The result is established by using some tools from the general
theory of processes such as Mertens decomposition of optional strong
(but not necessarily right continuous) supermartingales and some
tools from optimal stopping theory, as well as an appropriate
generalization of It\^{o}'s formula due to Gal'chouk and Lenglart. Two applications  on dynamic risk  measure and on optimal stopping will be  given.
\end{abstract}

\maketitle

\section{Introduction}

The notion of Backward Stochastic Differential Equations (BSDEs in
short) was introduced by \cite{Bismut(1973), Bismut(1976)} in the case of a linear driver. The nonlinear case
was developed by  \cite{PP1990, PP1992}. BSDEs have found
a number of applications in finance, that is pricing and hedging of
European options and recursive utilities (for
instance \cite{ELKarouiPeng1997}).\\
Reflected Backward Stochastic Differential Equations (RBSDEs in
short) have been introduced by
\cite{ELKaroui1997} and were useful, for example, in the study of
American option. The difference between the two types of equations
(BSDEs  and RBSDEs) is that the second can be seen as a variant of
the first in which the first component of the solution is
constrained to remain greater than or equal to a given process
called obstacle or barrier, and there is an additional nondecreasing
predictable process which keeps the first component of the solution
above the obstacle. The work of  \cite{ELKaroui1997} considers the case of a Brownian filtration and a continuous obstacle. After there have been several extensions
of this work to the case of a discontinuous obstacle, for example,
 \cite{Hamadene2002},
\cite{HamadeneOknine2003, HamadeneOknine2011},
\cite{Essaky2008} and \cite{Crepeymatoussi} ...\\
The right continuity of the obstacle is the difference between these
extensions and the paper of  \cite{Ouknine2015}. In this work, the
authors present a further extension of the theory of RBSDEs to the
case where the obstacle is not necessarily right-continuous in a
Brownian filtration.

In the present paper, we generalize the result of uniqueness and
existence of the RBSDE's solution in \cite{Ouknine2015} to the
case of a larger stochastic basis, i.e. in a  filtration that
supports a Brownian motion $W$ and an independent Poisson random
measure $\pi$, we establish existence and uniqueness of solutions,
in appropriate Banach spaces, to the following RBSDE:
\begin{multline}
Y_{\tau} =
\xi_{T}+\int_{\tau}^{T}f(s,Y_{s},Z_{s},\psi_{s})ds-\int_{\tau}^{T}Z_{s}dW_{s}-\int_{\tau}^{T}\int_{\mathcal{U}}\psi_{s}(u)\widetilde{\pi}(du,ds)
  -\int_{\tau}^{T}dM_{s}\\ +A_{T}-A_{\tau}+C_{T-}-C_{\tau-}\,\,\,\ \text{for all}\,\,\,\ \tau\in\mathcal{T}_{0,T}
   \label{mon0}
\end{multline}

The solution is given by $(Y,Z,\psi,M,A,C)$, where $M$  is an
orthogonal local martingale. We assume that the function $f$ is
Lipschitz  with respect to  $y$, $z$ and $\psi$. To prove our
results we use tools from the general theory of processes such as
Mertens decomposition of strong optional (but not necessarily
right-continuous) supermartingales (generalizing Doob-Meyer
decomposition) and some tools from optimal stopping theory, as well
as a generalization of It\^{o}'s formula to the case of strong optional
(but not necessarily right-continuous) semimartingales due to
\cite{Galchouk1981} and  \cite{Lenglart1980}.

We recover these natural differential equations by studying the limit of a system of reflected BSDEs

$$\left\{
        \begin{array}{ll}
          Y^{n}_{\tau}=\xi_{T}+\int_{\tau}^{T}f(s,Y^{n}_{s},Z^{n}_{s},\psi^{n}_{s})ds+K^{n}_{T}-K^{n}_{\tau} -\int_{\tau}^{T}Z^{n}_{s}dW_{s}-\int_{\tau}^{T}\int_{\mathcal{U}}\psi^{n}_{s}(u)\widetilde{\pi}(du,ds)\\
          -\int_{\tau}^{T}dM^{n}_{s} \\
          Y^{n}_{\tau} \geq \xi_{\tau}& \hbox{}
        \end{array}
      \right.
$$
Where $K^{n}_{t}=n\int_{0}^{t}(Y^{n}_{s}-\xi_{s})^{-}ds$. Essaky proved (in \cite{Essaky2008}), by a monotonic limit theorem, that $(Y^{n},Z^{n},K^{n},\psi^{n},M^{n})$ has, in some sense, a limit $(Y,Z,K,\psi,M)$  which satisfies a reflected BSDE with $\xi$  a c\`{a}dl\`{a}g barrier (see also \cite{Peng1999} for the case of filtration generated only by a brownian motion).

It is well known that if $\xi$ is a c\`{a}dl\`{a}g barrier then $Y$ is also a c\`{a}dl\`{a}g process  (Theorem 3.1 in  \cite{Essaky2008} for filtration generated by a Brownian motion
and Poisson point process, and Lemma 2.2  in \cite{Peng1999} for the Brownian filtration). But if the barrier $\xi$ is only optional the limit $Y$ of $Y^{n}$ is $\mathcal{E}^{f}$-super-martingale, then $Y$ has left and right limits (see \cite{DellacherieMeyer1980}, Theorem 4 page 408).

In this sense, we know that $(Y^{n},Z^{n},\psi^{n},M^{n})$ converge  to  $(Y,Z,\psi,M)$  and the limit $K$ of $K^{n}_{t}=n\int_{0}^{t}(Y^{n}_{s}-\xi_{s})^{-}ds$ is a l\`{a}dl\`{a}g process that can be written as $K=A+C_{-}$ where $A$ an increasing c\`{a}dl\`{a}g predictable process satisfying $A_{0}=0$, $E(A_{T})<\infty$, and $C$ an increasing c\`{a}dl\`{a}g  optional process and $E(C_{T})<\infty$.

The paper is decomposed as follows: in the second section, we give
the mathematical setting (preliminary, definitions and properties).
In subsection 2.1   we recall the change of variables formula for
optional semimartingales which are not necessarily right continuous
(Gal'chouk-Lenglart decomposition for strong optional
semimartingales). In the third  section, we define our RBSDE and we
prove existence and uniqueness of the solution in a general
filtration. In the last  section,  we give two applications of reflected BSDEs where the right-continuity of the obstacle is not necessarily used: application  on dynamic risk  measure and on optimal stopping.

\section{Preliminaries}
Let $T>0$ be a fixed positive real number. Let us consider a
filtered probability space $(\Omega,
\mathcal{F},\mathbb{P},\mathbb{F}=\{\mathcal{F}_{t}, t\geq 0\})$.
The filtration is assumed to be complete, right continuous and
quasi-left continuous, which means that for every sequence
($\tau_n$) of $\mathbb{F}$-stopping times such that $\tau_n\nearrow
\tau$ for some stopping time $\tau$ we have
$\bigvee_{n\in\N}\mathcal{F}_{\tau_n}=\mathcal{F}_{\tau}$. We assume
that $(\Omega, \mathcal{F},\mathbb{P},\mathbb{F}=\{\mathcal{F}_{t},
t\geq 0\})$ supports a $k$-dimensional Brownian motion $W$ and a
Poisson random measure $\pi$ with intensity  $\mu(du)dt$ on the
space $\mathcal{U}\subset\R^{m}\setminus\{0\}$. The measure $\mu$ is
$\sigma$-finite on $\mathcal{U}$ such that
\begin{equation}
\int_{\mathcal{U}}(1\wedge|u|^{2})\mu(du)<+\infty
 \label{moneq2}
\end{equation}
The compensated Poisson random measure $\pi$:
$\widetilde{\pi}(du,dt)=\pi(du,dt)-\mu(du)dt$ is a martingale w.r.t.
the filtration $\mathbb{F}$.\\
In this paper for a given $T>0$, we denote:
\begin{itemize}
\item $\mathcal{T}_{t,T}$ is the set of all  stopping times $\tau$
such that $\mathbb{P}(t\leq \tau \leq T)=1$. More generally, for a
given stopping time  $\nu$ in $\mathcal{T}_{0,T}$, we denote by
   $\mathcal{T}_{\nu,T}$ the set of all  stopping times $\tau$
  such that $\mathbb{P}(\nu \leq \tau \leq T)=1$.
\item  $\mathcal{P}$ is the predictable $\sigma$-field on $\Omega\times
[0,T]$ and
   $$\widetilde{\mathcal{P}}=\mathcal{P}\otimes\mathcal{B}(\mathcal{U})$$
   where $\mathcal{B}(\mathcal{U})$ is the Borelian $\sigma$-field on
   $\mathcal{U}$.
\item $L^{2}(\mathcal{F}_{T})$ is the set of random variables which are
$\mathcal{F}_{T}$-measurable and square-integrable.
\item On $\widetilde{\Omega}=\Omega\times[0,T]\times
    \mathcal{U}$, a function that is $\widetilde{\mathcal{P}}$-measurable, is called predictable.
\item $G_{loc}(\pi)$ is the set of $\widetilde{\mathcal{P}}$-measurable functions
         $\psi$ on $\widetilde{\Omega}$ such that for any  $t\geq 0$  a.s.
    $$\int_{0}^{t}\int_{\mathcal{U}}(|\psi_{s}(u)|^{2}\wedge|\psi_{s}(u)|)\mu(du)<+\infty.$$
\item $\mathds{H}^{2,T}$ is the set of real-valued predictable processes $\phi$ such that
         $$\|\phi\|_{\mathds{H}^{2,T}}^{2}=E\Bigl[\int_{0}^{T}|\phi_{t}|^{2}dt\Bigr]<\infty.$$

\item $\mathcal{M}_{loc}$ is the set of c\`{a}dl\`{a}g local martingales orthogonal to $W$ and $\widetilde{\pi}$:
 if $M\in\mathcal{M}_{loc}$ then  $$[M,W^{i}]_{t}=0, \,\,\,\  1\leq i\leq k,
   \,\,\,\,\ [M,\widetilde{\pi}(A)]_{t}=0$$ for all
   $A\in\mathcal{B}(\mathcal{U})$.
\item $\mathcal{M}$ is the subspace of $\mathcal{M}_{loc}$ of martingales.
\end{itemize}
As explained above, the filtration  $\mathbb{F}$ supports the
Brownian motion $W$ and the Poisson random measure $\pi$. We have
the following lemma that we  can  find in
\cite{JacodShiryaev2003} (Chapter III, Lemma 4.24):

\begin{lemma}
$\label{lemme1}$ Every local martingale $N$ has a decomposition
\begin{equation}
N_{t}=\int_{0}^{t}Z_{s}dW_{s}+\int_{0}^{t}\int_{\mathcal{U}}\psi_{s}(u)\widetilde{\pi}(du,ds)+M_{t}
 \label{moneq3}
\end{equation}
where  $M\in\mathcal{M}_{loc}$, $Z\in\mathds{H}^{2,T}$ and $\psi\in
G_{loc}(\mu).$
\end{lemma}
Now to define the solution of our reflected backward stochastic
differential equation (RBSDE),  let us introduce the following
spaces:
\begin{itemize}
\item  $\mathds{S}^{2,T}$ is the set of real-valued optional processes $\phi$ such
that:
         $$\||\phi\||_{\mathds{S}^{2,T}}^{2}=E\Bigl[ess\sup_{\tau\in\mathcal{T}_{0,T}}|\phi_{\tau}|^{2}\Bigr]<\infty.$$
\item $\mathbb{M}^{2}$ is the subspace of $\mathcal{M}$ of all martingales such that:
         $$\|M\|_{\mathbb{M}^{2}}^{2}=E([M,M]_{T})=E([M]_{T})<+\infty.$$
\item $\mathbb{L}^{2}_{\pi}(0,T)=\mathbb{L}^{2}_{\pi}(\Omega\times(0,T)\times\mathcal{U})$
         is the set of all processes $\psi\in G_{loc}(\mu)$ such that:
         $$\|\psi\|^{2}_{\mathbb{L}^{2}_{\pi}}=E\Bigl[\int_{0}^{T}\int_{\mathcal{U}}|\psi_{s}(u)|^{2}\mu(du)ds\Bigr]<+\infty$$
\item  $\mathbb{L}^{2}_{\mu}(0,T)=\mathbb{L}^{2}(\mathcal{U},\mu;\R^{d})$
        is the set of all measurable functions $\psi:\mathcal{U}\longrightarrow
         \R^{d}$ such that:
          $$\|\psi\|_{\mathbb{L}^{p}_{\mu}}^{2}=\int_{\mathcal{U}}|\psi(u)|^{2}\mu(du)<+\infty.$$
\item $\mathcal{E}^{2}(0,T)=\mathds{S}^{2,T}\times\mathds{H}^{2,T}\times\mathbb{L}^{2}_{\pi}(0,T)
\times\mathbb{M}^{2}\times\mathds{S}^{2,T}\times\mathds{S}^{2,T}.$
\end{itemize}

The random variable $\xi$ is $\mathcal{F}_{T}$-measurable with
values in $\R^{d}$ $(d\geq 1)$ and
$f:\Omega\times[0,T]\times\R^{d}\times\R^{d\times
k}\times\mathbb{L}^{2}_{\mu}\longrightarrow\R^{d}$ is a random
function measurable with respect to
$Prog\times\mathcal{B}(\R^{d})\times\mathcal{B}(\R^{d\times
k})\times\mathcal{B}(\mathbb{L}^{2}_{\mu})$ where $Prog$ denotes the
$\sigma$-field of progressive subsets
of $\Omega\times[0,T]$.\\
In the following we  denote the spaces $\mathds{H}^{2,T}$ and
$\mathds{S}^{2,T}$ by $\mathds{H}^{2}$ and $\mathds{S}^{2}$, as well
as the norms $\parallel.\parallel_{\mathds{H}^{2,T}}$ and
$\parallel\mid.\parallel\mid_{\mathds{S}^{2,T}}$ by
$\parallel.\parallel_{\mathds{H}^{2}}$ and
$\parallel\mid.\parallel\mid_{\mathds{S}^{2}}$.
\begin{definition}
A function $f$ is said to be a driver if:
\begin{itemize}
    \item $f:\Omega\times [0,T]\times \R^{2}\times\mathbb{L}^{2}_{\mu}\longrightarrow \R$ \\
    $(\omega,t,y,z,\psi)\longmapsto f(\omega,t,y,z,\psi)$ is $\mathcal{P}\otimes\mathcal{B}(\mathds{R}^{2})\otimes\mathcal{B}(\mathbb{L}^{2}_{\mu})$-measurable.
    \item $E[\int_{0}^{T}\mid f(t,0,0,0)\mid^{2}dt]<\infty$.
\end{itemize}
A driver $f$ is called a Lipschitz driver if moreover there exists a
constant $K\geq0$ such that $\mathbb{P}\otimes dt$-a.s., for each
$(y_1,z_1,\psi_1)$ and $(y_2,z_2,\psi_2)$
\begin{equation}
\Bigl|f(\omega,t,y_1,z_1,\psi_1)-f(\omega,t,y_2,z_2,\psi_2)\Bigr|\leq
K\Bigl(|y_1-y_2|+|z_1-z_2|+\|\psi_1-\psi_2\|_{\mathbb{L}^{2}_{\mu}}\Bigr)
 \label{moneq0}
\end{equation}
\end{definition}
For a l\`{a}dl\`{a}g process $\phi$, we denote by $\phi_{t+}$ and
$\phi_{t-}$ the right-hand and left-hand limit of $\phi$ at  $t$. We
denote by $\Delta_{+}\phi_{t}=\phi_{t+}-\phi_{t}$ the size of the
right jump of $\phi$ at $t$, and by
$\Delta\phi_{t}=\phi_{t}-\phi_{t-}$ the size of the left jump of
$\phi$ at $t$.

 We give a useful property of the space
$\mathds{S}^{2}$:
\begin{proposition}
The space $\mathds{S}^{2}$ endowed with the norm
$\parallel\mid.\parallel\mid_{\mathds{S}^{2}}$ is a Banach space.
\end{proposition}
\begin{proof}
The proof is given in  \cite{Ouknine2015} (Proposition 2.1).
\end{proof}
The following proposition can be found in
\cite{AshkanNikeghbali2006} (Theorem 3.2.).
\begin{proposition}
Let $(X_t)$ and $(Y_t)$ be two optional processes. If for every
finite stopping time $\tau$ one has, $X_{\tau}=Y_{\tau}$, then the
processes $(X_t)$ and $(Y_t)$ are indistinguishable.
\end{proposition}

Let $\beta>0$. We will also use the following notation: \\
 For  $\phi\in\mathds{H}^{2}$,
$\|\phi\|_{\beta}^{2}:=E[\int_{0}^{T}e^{\beta s}\phi^{2}_{s}ds]$. We
note that on the space $\mathds{H}^{2}$ the norms $\|.\|_{\beta}$
and  $\|.\|_{\mathds{H}^{2}}$ are equivalent.\\ For
$\phi\in\mathds{S}^{2}$, we define
$\||\phi\||_{\beta}^{2}:=E[ess\sup_{\tau\in\mathcal{T}_{0,T}}e^{\beta
\tau}|\phi_{\tau}|^{2}]$. We note that $\||.|\|_{\beta}$ is a norm
on $\mathds{S}^{2}$ equivalent to the norm $\||.|\|_{\mathds{S}^{2}}$.\\
 For $\phi\in\mathbb{L}^{2}_{\pi}$, the defined norm
 $\|\phi\|_{\mathbb{L}^{2,\beta}_{\pi}}=\sqrt{E[\int_{0}^{T}e^{\beta
 s}\int_{\mathcal{U}}|\phi_{s}(u)|^{2}\mu(du)ds]}$ is equivalent to the norm
 $\|\phi\|_{\mathbb{L}^{2}_{\pi}}$ on
 $\mathbb{L}^{2}_{\pi}$.\\
For  $M\in\mathbb{M}^{2}$, we have the equivalence between
 $\|M\|_{\mathbb{M}_{\beta}^{2}}=\sqrt{E[\int_{0}^{T}e^{\beta
 s}d[M]_{s}]}$ and
 $\|M\|_{\mathbb{M}^{2}}$ on
 $\mathbb{M}^{2}$.
\subsection{Gal'chouk-Lenglart decomposition for strong optional semimartingales.}
In this section, we recall the change of variables formula for
optional semimartingales which are not necessarily cad. The result
can be seen as a generalization of the classical It\^{o} formula and can
be found in (\cite{Galchouk1981}, (Theorem 8.2)),
(\cite{Lenglart1980},(Section 3, page 538)). We recall the result
in our framework in which the underlying filtered probability space
satisfies the usual conditions.
\begin{theorem}
\label{galchouklenglart} \textbf{(Gal'chouk-Lenglart)}  Let
$n\in\N$. Let $X$ be an n-dimensional optional semimartingale, i.e.
$X = (X_1,... ,X_n)$ is an n-dimensional optional process with
decomposition $X_{t}^{k}=X_{0}^{k}+N_{t}^{k}+A_{t}^{k}+B_{t}^{k}$,
for all $k\in\{1,...,n\}$ where $N_{t}^{k}$ is a (c\`{a}dl\`{a}g) local
martingale, $A_{t}^{k}$ is a right-continuous process of finite
variation such that $A_{0}=0$ and $B_{t}^{k}$ is a left-continuous
process of finite variation which is purely discontinuous and such
that $B_{0-}=0$. Let $F$ be a twice continuously differentiable
function on $\R^{n}$. Then, almost surely, for all $t\geq0$,
\begin{eqnarray*}
  F(X_{t}) &=& F(X_{0})+\sum_{k=1}^{n}\int_{]0,t]}D^{k}F(X_{s-})d(N^{k}+A^{k})_{s} \\
           &+& \frac{1}{2}\sum_{k,l=1}^{n}\int_{]0,t]}D^{k}D^{l}F(X_{s-})d<N^{kc},N^{lc}>_{s}\\
           &+& \sum_{0<s\leq t}\Bigl[F(X_{s})-F(X_{s-})-\sum_{k=1}^{n}D^{k}F(X_{s-})\Delta X^{k}_{s}\Bigr]\\
           &+& \sum_{k=1}^{n}\int_{[0,t[}D^{k}F(X_{s})d(B^{k})_{s+} \\
           &+& \sum_{0\leq s<t}\Bigl[F(X_{s+})-F(X_{s})-\sum_{k=1}^{n}D^{k}F(X_{s})\Delta_{+}X^{k}_{s}\Bigr]
\end{eqnarray*}
where $D^k$ denotes the differentiation operator with respect to the
$k$-th coordinate, and $N^{kc}$ denotes the continuous part of
$N^k$.
\end{theorem}
\begin{corollary}
\label{corollaire1} Let $Y$ be a one-dimensional optional
semimartingale with decomposition $Y_t=Y_0+N_t+A_t+B_t$, where $N$,
$A$ and $B$ are as in the above theorem. Let $\beta
>0$. Then, almost surely, for all $t$ in $[0,T]$,
\begin{eqnarray*}
  e^{\beta t}Y_{t}^{2} &=&  Y_{0}^{2}+\int_{0}^{t}\beta e^{\beta s}Y_{s}^{2}ds+2\int_{0}^{t}e^{\beta s}Y_{s-}d(A+N)_{s}\\
                       &+& \int_{0}^{t} e^{\beta s}d<N^{c},N^{c}>_{s}\\
                       &+& \sum_{0<s\leq t}e^{\beta s}(Y_{s}-Y_{s-})^{2}+2\int_{0}^{t}e^{\beta
s}Y_{s}d(B)_{s+}+\sum_{0\leq s<t}e^{\beta
  s}(Y_{s+}-Y_{s})^{2}
\end{eqnarray*}
\end{corollary}
\begin{proof}
For the corollary demonstration, it suffices to apply the change of
variables formula from Theorem$~\ref{galchouklenglart}$ with $n=2$,
$F(x,y)=xy^{2}$, $X_{t}^{1}= e^{\beta t}$ and  $X_{t}^{2}=Y_{t}$.
Indeed, by applying Theorem$~\ref{galchouklenglart}$ and by noting
that the local martingale part and the purely discontinuous part of
$X^{1}$ are both equal to $0$, we obtain
\begin{eqnarray*}
  e^{\beta t}Y_{t}^{2} &=&  Y_{0}^{2}+\int_{0}^{t}\beta e^{\beta s}Y_{s}^{2}ds+2\int_{0}^{t}e^{\beta s}Y_{s-}d(A+N)_{s}\\
                       &+& \int_{0}^{t} e^{\beta s}d<N^{c},N^{c}>_{s}\\
                       &+& \sum_{0<s\leq t}e^{\beta s}(Y_{s}^{2}-Y_{s-}^{2}-2Y_{s-}(Y_{s}-Y_{s-}))+2\int_{0}^{t}e^{\beta s}Y_{s}d(B)_{s+}\\
                       &+&  \sum_{0\leq s<t}e^{\beta s}(Y_{s+}^{2}-Y_{s}^{2}-2Y_{s}(Y_{s+}-Y_{s}))
\end{eqnarray*}
The desired expression follows as
$(Y_{s}-Y_{s-})^{2}=Y_{s}^{2}-Y_{s-}^{2}-2Y_{s-}(Y_{s}-Y_{s-})$ and
$(Y_{s+}-Y_{s})^{2}=Y_{s+}^{2}-Y_{s}^{2}-2Y_{s}(Y_{s+}-Y_{s})$.
\end{proof}

\section{RBSDEs whose obstacles are not c\`{a}dl\`{a}g in a general filtration.}

Let $T>0$ be a fixed terminal time. Let $f$ be a driver. Let
$\xi=(\xi_t)_{t\in[0,T]}$ be a left-limited process in
$\mathds{S}^{2}$. We suppose moreover that the process $\xi$ is
right upper-semicontinuous (r.u.s.c. for short). A process $\xi$
satisfying the previous properties will be called a barrier, or an
obstacle.
\begin{definition}
\label{definition1} A process $(Y,Z,\psi,M,A,C)$ is said to be a
solution to the reflected BSDE with parameters $(f,\xi)$, where $f$
is a driver and $\xi$ is an obstacle, if $(Y,Z,\psi,M,A,C)
\in\mathcal{E}^{2}(0,T)$ and

\begin{multline}
Y_{\tau} =
\xi_{T}+\int_{\tau}^{T}f(s,Y_{s},Z_{s},\psi_{s})ds-\int_{\tau}^{T}Z_{s}dW_{s}-\int_{\tau}^{T}\int_{\mathcal{U}}\psi_{s}(u)\widetilde{\pi}(du,ds)
  -\int_{\tau}^{T}dM_{s}\\ +A_{T}-A_{\tau}+C_{T-}-C_{\tau-} \,\,\,\ \text{for all}\,\,\,\ \tau\in\mathcal{T}_{0,T}
   \label{moneq4}
\end{multline}
\begin{equation}
 Y\geq \xi \,\,\ \text{(up to an evanescent set)} \,\,\
\text{a.s.}
 \label{moneq5}
\end{equation}
\begin{equation}
M\in\mathcal{M}_{loc} \,\,\ \text{and} \,\,\ M_0=0.
 \label{martingale}
\end{equation}
In the above, the process $A$ is a nondecreasing right-continuous predictable
process with $A_0=0$, $E(A_T)<\infty$  such that:
\begin{equation}
\int_{0}^{T}1_{ \{Y_{t}>\xi_{t}\}}dA_{t}^{c}=0 \,\,\  a.s.\,\,\ and\,\,\ (Y_{\tau-}-
\xi_{\tau-})(A_{\tau}^{d}-A_{\tau-}^{d})=0 \,\,\  a.s. \,\,\ \forall
(predictable)  \tau\in\mathcal{T}_{0,T}
 \label{moneq6}
\end{equation}
And the process $C$ is a nondecreasing right-continuous adapted purely
discontinuous process with  $C_{0-}=0$, $E(C_T)<\infty$
such that:
\begin{equation}
 (Y_{\tau}-\xi_{\tau})(C_{\tau}-C_{\tau-})=0 \,\,\
a.s.\,\,\ \forall \tau\in\mathcal{T}_{0,T}
 \label{moneq7}
\end{equation}
Here $A^{c}$  denotes the continuous part of the nondecreasing
process $A$ and $A^{d}$ its discontinuous part.
\end{definition}
\begin{remark}
We note that a process $(Y,Z,\psi,M,A,C) \in\mathcal{E}^{2}(0,T)$
satisfies equation $~\eqref{moneq4}$ in the above definition if and
only if  $\forall  t\in[0,T]$, a.s.
\begin{eqnarray*}
  Y_{t} = \xi_{T}+\int_{t}^{T}f(s,Y_{s},Z_{s},\psi_{s})ds-\int_{t}^{T}Z_{s}dW_{s}&-&\int_{t}^{T}\int_{\mathcal{U}}\psi_{s}(u)\widetilde{\pi}(du,ds)
- \int_{t}^{T}dM_{s} \\
           &+& A_{T}-A_{t}+C_{T-}-C_{t-}
\end{eqnarray*}
\end{remark}
\begin{remark} If $(Y,Z,\psi,M,A,C) \in\mathcal{E}^{2}(0,T)$ satisfies
the above definition, then the process $Y$ has left and right
limits. Moreover, the process
$(Y_{t}+\int_{0}^{t}f(s,Y_{s},Z_{s},\psi_{s})ds)_{t\in[0,T]}$ is a
strong supermartingale.
\end{remark}
The proof of the existence and uniqueness of the reflected BSDE
solution defined above is based on a useful result (following lemma)
in the case of $f$ depends only on $s$ and $\omega$ (i.e.
$f(s,y,z,\psi)=f(s,\omega)$), the corollary$~\ref{corollaire1}$ and
the lemma$~\ref{lemme1}$. To this purpose, we first prove a lemma
which will be used in the sequel.
\begin{lemma}
\label{lemme2} Let
$(Y^{1},Z^{1},\psi^{1},M^{1},A^{1},C^{1})\in\mathcal{E}^{2}(0,T)$
(resp.
$(Y^{2},Z^{2},\psi^{2},M^{2},A^{2},C^{2})\in\mathcal{E}^{2}(0,T)$.)
be a solution to the RBSDE associated with driver $f^{1}(s,\omega)$
(resp.$f^{2}(s,\omega)$) and with obstacle $\xi$. There exists $c>0$
such that for all $\epsilon>0$, for all
$\beta>\frac{1}{\epsilon^{2}}$ we have
\begin{equation}
\|Z^{1}-Z^{2}\|^{2}_{\beta}+\|M^{1}-M^{2}\|_{\mathbb{M}_{\beta}^{2}}^{2}
    +\|\psi^{1}-\psi^{2}\|^{2}_{\mathbb{L}^{2,\beta}_{\pi}}\leq
\epsilon^{2}\|f^{1}-f^{2}\|_{\beta}^{2}
 \label{unicite}
 \end{equation}
 and
\begin{equation}
\||Y^{1}-Y^{2}|\|_{\beta}^{2}\leq
4\epsilon^{2}(1+4c^{2})\|f^{1}-f^{2}\|_{\beta}^{2} \label{moneq8}
\end{equation}
\end{lemma}
\begin{proof}
Let  $\beta >0$ and   $\epsilon>0$ be such that $\beta\geq
\frac{1}{\epsilon^{2}}$. We set  $\widetilde{Y}:=Y^{1}-Y^{2}$,
$\widetilde{Z}:=Z^{1}-Z^{2}$, $\widetilde{\psi}:=\psi^{1}-\psi^{2}$,
$\widetilde{M}:=M^{1}-M^{2}$, $\widetilde{A}:=A^{1}-A^{2}$,
$\widetilde{C}:=C^{1}-C^{2}$ and
$\widetilde{f}(\omega,t):=f^{1}(\omega,t)-f^{2}(\omega,t)$. We note
that  $\widetilde{Y}_{T}:=\xi_{T}-\xi_{T}=0$. Moreover,
\begin{multline}
\widetilde{Y}_{\tau}=\int_{\tau}^{T}\widetilde{f}(s)ds-\int_{\tau}^{T}\widetilde{Z}_{s}dW_{s}
-\int_{\tau}^{T}\int_{\mathcal{U}}\widetilde{\psi}_{s}(u)\widetilde{\pi}(du,ds)-\widetilde{M}_{T}+\widetilde{M}_{\tau}
+\widetilde{A}_{T}-\widetilde{A}_{\tau}\\+\widetilde{C}_{T-}-\widetilde{C}_{\tau-}
\,\,\ a.s.\,\,\ \forall \tau\in\mathcal{T}_{0,T}, \label{moneq9}
\end{multline}
i.e.
$$\widetilde{Y}_{\tau}=\widetilde{Y}_{0}-\int_{0}^{\tau}\widetilde{f}(s)ds+
\int_{0}^{\tau}\widetilde{Z}_{s}dW_{s}
+\int_{0}^{\tau}\int_{\mathcal{U}}\widetilde{\psi}_{s}(u)\widetilde{\pi}(du,ds)+\widetilde{M}_{\tau}-\widetilde{A}_{\tau}-\widetilde{C}_{\tau-}\,\,\
a.s.\,\,\ \forall \tau\in\mathcal{T}_{0,T},$$

Since $\widetilde{M}_{0}=M_{0}^{1}-M_{0}^{2}=0$,
$\widetilde{A}_{0}=A_{0}^{1}-A_{0}^{2}=0$ and
$\widetilde{C}_{0-}=C_{0-}^{1}-C_{0-}^{2}=0$. Thus we see that
$\widetilde{Y}$  is an optional (strong) semimartingale with
decomposition
$\widetilde{Y}_{t}=\widetilde{Y}_{0}+N_{t}+A_{t}+B_{t}$, where
$N_{t}=\int_{0}^{t}\widetilde{Z}_{s}dW_{s}+\int_{0}^{t}\int_{\mathcal{U}}\widetilde{\psi}_{s}(u)\widetilde{\pi}(du,ds)+\widetilde{M}_{t}$,
$A_{t}=-\int_{0}^{t}\widetilde{f}(s)ds-\widetilde{A}_{t}$ and
$B_{t}=-\widetilde{C}_{t-}$ (the notation is that of
(\ref{galchouklenglart})), Applying Corollary$~\ref{corollaire1}$
to $\widetilde{Y}$ gives: almost surely, for all $t\in[0,T]$,
\begin{eqnarray*}
  e^{\beta t}\widetilde{Y}_{t}^{2} &=&  -\int_{0}^{t}\beta e^{\beta s}\widetilde{Y}_{s}^{2}ds+2\int_{0}^{t}e^{\beta s}\widetilde{Y}_{s-}d(A+N)_{s}\\
   &-& \int_{0}^{t} e^{\beta s}d<N^{c},N^{c}>_{s} \\
   &-& \sum_{0< s\leq t}e^{\beta s}(\widetilde{Y}_{s}-\widetilde{Y}_{s-})^{2}-\int_{0}^{t}2e^{\beta
s}\widetilde{Y}_{s}d(B)_{s+}-\sum_{0\leq s<t}e^{\beta
  s}(\widetilde{Y}_{s+}-\widetilde{Y}_{s})^{2}
\end{eqnarray*}
Using the expressions of $N$, $A$ and $B$  and the fact that
$\widetilde{Y}_{T}=0$, we get: almost surely, for all $t\in[0,T]$,
\begin{eqnarray*}
   e^{\beta t}\widetilde{Y}_{t}^{2}+\int_{t}^{T} e^{\beta s}d<N^{c},N^{c}>_{s} &=&
   -\int_{t}^{T}\beta e^{\beta s}\widetilde{Y}_{s}^{2}ds+2\int_{t}^{T}e^{\beta
   s}\widetilde{Y}_{s-}\widetilde{f}(s)ds \\
  &+& 2\int_{t}^{T}e^{\beta s}\widetilde{Y}_{s-}d\widetilde{A}- 2\int_{t}^{T}e^{\beta s}\widetilde{Y}_{s-}\widetilde{Z}_{s}dW_{s}\\
  &-&  2\int_{t}^{T}e^{\beta s}\widetilde{Y}_{s-}\int_{\mathcal{U}}\widetilde{\psi}_{s}(u)\widetilde{\pi}(du,ds)- 2\int_{t}^{T}e^{\beta s}\widetilde{Y}_{s-}d\widetilde{M}_{s} \\
  &-& \sum_{t<s\leq T}e^{\beta
s}(\widetilde{Y}_{s}-\widetilde{Y}_{s-})^{2}+\int_{t}^{T}2e^{\beta
s}\widetilde{Y}_{s}d(\widetilde{C})_{s}\\
&-&\sum_{t\leq s<T}e^{\beta
  s}(\widetilde{Y}_{s+}-\widetilde{Y}_{s})^{2}
\end{eqnarray*}
Then
\begin{eqnarray*}
   e^{\beta t}\widetilde{Y}_{t}^{2} + \int_{t}^{T} e^{\beta s}\widetilde{Z}^{2}_{s}ds &+& \int_{t}^{T} e^{\beta s}\int_{\mathcal{U}}|\widetilde{\psi}_{s}(u)|^{2}\mu(du)ds+\int_{t}^{T} e^{\beta
   s}d<\widetilde{M}^{c},\widetilde{M}^{c}>_{s}=\\
   &-&\int_{t}^{T}\beta e^{\beta s}\widetilde{Y}_{s}^{2}ds+2\int_{t}^{T}e^{\beta s}\widetilde{Y}_{s-}\widetilde{f}(s)ds \\
  &+& 2\int_{t}^{T}e^{\beta s}\widetilde{Y}_{s-}d\widetilde{A}- 2\int_{t}^{T}e^{\beta s}\widetilde{Y}_{s-}\widetilde{Z}_{s}dW_{s}\\
  &-&  2\int_{t}^{T}e^{\beta s}\widetilde{Y}_{s-}\int_{\mathcal{U}}\widetilde{\psi}_{s}(u)\widetilde{\pi}(du,ds)- 2\int_{t}^{T}e^{\beta s}\widetilde{Y}_{s-}d\widetilde{M}_{s} \\
  &-& \sum_{t<s\leq T}e^{\beta
s}(\widetilde{Y}_{s}-\widetilde{Y}_{s-})^{2}+\int_{t}^{T}2e^{\beta
s}\widetilde{Y}_{s}d(\widetilde{C})_{s}\\
&-&\sum_{t\leq s<T}e^{\beta s}(\widetilde{Y}_{s+}-\widetilde{Y}_{s})^{2}
\end{eqnarray*}
It is clear that for all $t\in[0,T]$  $-\sum_{t<s\leq T}e^{\beta
s}(\widetilde{Y}_{s}-\widetilde{Y}_{s-})^{2}-\sum_{t\leq
s<T}e^{\beta s}(\widetilde{Y}_{s+}-\widetilde{Y}_{s})^{2} \leq 0$.
By applying the inequality $2ab\leq
(\frac{a}{\epsilon})^{2}+\epsilon^{2}b^{2}$, valid for all $(a,b)$
in $\mathbb{R}^{2}$, we get: a.e. for all $t\in[0,T]$
\begin{eqnarray*}
  -\int_{t}^{T}\beta e^{\beta s}\widetilde{Y}_{s}^{2}ds+2\int_{t}^{T}e^{\beta
   s}\widetilde{Y}_{s-}\widetilde{f}(s)ds &\leq & -\int_{t}^{T}\beta e^{\beta
   s}\widetilde{Y}_{s}^{2}ds +\frac{1}{\epsilon^{2}}\int_{t}^{T}e^{\beta
   s}\widetilde{Y}_{s-}^{2}ds\\
   &+&\epsilon^{2}\int_{t}^{T}e^{\beta s}\widetilde{f}(s)^{2}ds\\
   &=& (\frac{1}{\epsilon^{2}}-\beta)\int_{t}^{T}e^{\beta
   s}\widetilde{Y}_{s-}^{2}ds+\epsilon^{2}\int_{t}^{T}e^{\beta
   s}\widetilde{f}(s)^{2}ds
\end{eqnarray*}
As $\beta\geq\frac{1}{\epsilon^{2}}$, we have
$(\frac{1}{\epsilon^{2}}-\beta)\int_{t}^{T}e^{\beta
   s}\widetilde{Y}_{s-}^{2}ds\leq0$ for all $t\in[0,T]$ a.s.\\
  Next, we  have also that the term  $\int_{t}^{T}e^{\beta s}\widetilde{Y}_{s}d\widetilde{C}_{s}$
   is non-positive. Indeed a.s. for all $t\in[0,T]$,
   $$\int_{t}^{T}e^{\beta s}\widetilde{Y}_{s}d\widetilde{C}_{s}
   =\sum_{t\leq s<T}e^{\beta
   s}\widetilde{Y}_{s}\bigtriangleup\widetilde{C}_{s}$$ and a.s. for all $t\in[0,T]$
\begin{equation}
\widetilde{Y}_{t}\bigtriangleup\widetilde{C}_{t}=(Y^{1}_{t}-Y^{2}_{t})\bigtriangleup
C^{1}_{t}-(Y^{1}_{t}-Y^{2}_{t})\bigtriangleup C^{2}_{t}
 \label{moneq10}
\end{equation}
We use property $~\eqref{moneq7}$ of $C^1$ and the fact that
$Y^2\geq \xi$ to obtain: a.s. for all $t\in[0,T]$
$$(Y^{1}_{t}-Y^{2}_{t})\bigtriangleup
C^{1}_{t}=(Y^{1}_{t}-\xi_{t})\bigtriangleup
C^{1}_{t}-(Y^{2}_{t}-\xi_{t})\bigtriangleup
C^{1}_{t}=0-(Y^{2}_{t}-\xi_{t})\bigtriangleup C^{1}_{t}\leq0$$
Similarly, we obtain: a.s. for all $t\in[0,T]$,
$$(Y^{1}_{t}-Y^{2}_{t})\bigtriangleup
C^{2}_{t}=(Y^{1}_{t}-\xi_{t})\bigtriangleup
C^{2}_{t}-(Y^{2}_{t}-\xi_{t})\bigtriangleup
C^{2}_{t}=(Y^{1}_{t}-\xi_{t})\bigtriangleup C^{2}_{t}-0\geq0.$$ We
also show that $\int_{t}^{T}e^{\beta
s}\widetilde{Y}_{s-}d\widetilde{A}$ is non-positive by using
property $~\eqref{moneq6}$ of the definition of the RBSDE and the
fact that $Y^i\geq \xi$ for $i=1,2$ and that $A^i=A^{i,c}+A^{i,d}$
(see also \cite{QuenezSulem2014}). Then
\begin{multline}
   e^{\beta t}\widetilde{Y}_{t}^{2}+\int_{t}^{T} e^{\beta s}\widetilde{Z}^{2}_{s}ds+\int_{t}^{T} e^{\beta s}d<\widetilde{M}^{c},\widetilde{M}^{c}>_{s}
    +\int_{t}^{T} e^{\beta s}\int_{\mathcal{U}}|\widetilde{\psi}_{s}(u)|^{2}\mu(du)ds\leq \\
    \epsilon^{2}\int_{t}^{T}e^{\beta s}\widetilde{f}^{2}(s)ds-2\int_{t}^{T}e^{\beta s}\widetilde{Y}_{s-}\widetilde{Z}_{s}dW_{s}\\
     -2\int_{t}^{T}\int_{\mathcal{U}}e^{\beta s}\widetilde{Y}_{s-}\widetilde{\psi}_{s}(u)\widetilde{\pi}(du,ds)-
     2\int_{t}^{T}e^{\beta s}\widetilde{Y}_{s-}d\widetilde{M}_{s}
 \,\,\ \forall \,\,\ a.s. \,\,\ t\in[0,T].
 \label{moneq11}
\end{multline}
We now show that the term $\int_{0}^{T}e^{\beta
s}\widetilde{Y}_{s-}\widetilde{Z}_{s}dW_{s}$ has zero expectation.
To this purpose, we show that $E[\sqrt{\int_{0}^{T}e^{2\beta
s}\widetilde{Y}_{s-}^{2}\widetilde{Z}^{2}_{s}ds}]<\infty$, in the
same way that in the proof of Lemma 3.2 (A priori estimates) in
\cite{Ouknine2015}. By using the left-continuity of a.e.
trajectory of the process  $(\widetilde{Y}_{s-})$, we have
\begin{equation}
(\widetilde{Y}_{s-})^{2}(\omega)\leq
\sup_{t\in\mathbb{Q}}(\widetilde{Y}_{t-})^{2}(\omega) \,\,\
\text{for all} \,\,\ s\in(0,T],\,\,\  \text{for a.s.}\,\,\
\omega\in\Omega
 \label{moneq12}
\end{equation}
On the other hand, for all $t\in(0,T]$, a.s.,
$(\widetilde{Y}_{t-})^{2} \leq ess\sup_{\tau\in\mathcal{T}_{0,T}}(\widetilde{Y}_{\tau})^{2}$. Then
\begin{equation}
\sup_{t\in\mathbb{Q}}(\widetilde{Y}_{t-})^{2}\leq
ess\sup_{\tau\in\mathcal{T}_{0,T}}(\widetilde{Y}_{\tau})^{2} \,\,\
a.s.
 \label{moneq13}
\end{equation}
According to $~\eqref{moneq12}$ and $~\eqref{moneq13}$ we obtain
\begin{equation}
\int_{0}^{T}e^{2\beta
s}\widetilde{Y}_{s-}^{2}\widetilde{Z}_{s}^{2}ds\leq
\int_{0}^{T}e^{2\beta
s}\sup_{t\in\mathbb{Q}}\widetilde{Y}_{t-}^{2}\widetilde{Z}_{s}^{2}ds
\leq\int_{0}^{T}e^{2\beta
s}ess\sup_{\tau\in\mathcal{T}_{0,T}}\widetilde{Y}_{\tau}^{2}\widetilde{Z}_{s}^{2}ds
 \label{moneq14}
\end{equation}
Using $~\eqref{moneq14}$, together with Cauchy-Schwarz inequality,
gives
$$E\Bigl[\sqrt{\int_{0}^{T}e^{2\beta
s}\widetilde{Y}_{s-}^{2}\widetilde{Z}_{s}^{2}ds}\Bigr] \leq
E\Bigl[\sqrt{ess\sup_{\tau\in\mathcal{T}_{0,T}}\widetilde{Y}_{\tau}^{2}}\sqrt{\int_{0}^{T}e^{2\beta
s}\widetilde{Z}_{s}^{2}ds}\Bigr]$$ Then
\begin{equation}
  E\Bigl[\sqrt{\int_{0}^{T}e^{2\beta
s}\widetilde{Y}_{s-}^{2}\widetilde{Z}_{s}^{2}ds}\Bigr] \leq
\||\widetilde{Y}\||_{\mathds{S}^{2}}\|\widetilde{Z}\|_{2\beta}.
   \label{moneq15}
\end{equation}
We conclude that $ E\Bigl[\sqrt{\int_{0}^{T}e^{2\beta
s}\widetilde{Y}_{s-}^{2}\widetilde{Z}_{s}^{2}ds}\Bigr]<\infty$,
whence,  we get $E[\int_{0}^{T}e^{\beta
s}\widetilde{Y}_{s-}\widetilde{Z}_{s}dW_{s}]=0$. Next we show that
$E\Bigl[\int_{0}^{T}\int_{\mathcal{U}}e^{\beta
s}\widetilde{Y}_{s-}\widetilde{\psi}_{s}(u)\widetilde{\pi}(du,ds)\Bigr]=0$.
For this purpose, we first prove that
$E[\sqrt{\int_{0}^{T}\int_{\mathcal{U}}e^{2\beta
s}\widetilde{Y}_{s-}^{2}\widetilde{\psi}^{2}_{s}(u)\mu(du)ds}]<\infty$.
According to $~\eqref{moneq12}$ and $~\eqref{moneq13}$, we have
\begin{multline}
\int_{0}^{T}\int_{\mathcal{U}}e^{2\beta
s}\widetilde{Y}_{s-}^{2}\widetilde{\psi}^{2}_{s}(u)\mu(du)ds\leq
\int_{0}^{T}\int_{\mathcal{U}}e^{2\beta
s}\sup_{t\in\mathbb{Q}}\widetilde{Y}_{t-}^{2}\widetilde{\psi}^{2}_{s}(u)\mu(du)ds\\
\leq \int_{0}^{T}\int_{\mathcal{U}}e^{2\beta
s}ess\sup_{\tau\in\mathcal{T}_{0,T}}\widetilde{Y}_{\tau}^{2}\widetilde{\psi}^{2}_{s}(u)\mu(du)ds
 \label{moneq16}
\end{multline}
Using $~\eqref{moneq16}$ and  Cauchy-Schwarz inequality, gives
$$E\Bigl[\sqrt{\int_{0}^{T}\int_{\mathcal{U}}e^{2\beta
s}\widetilde{Y}_{s-}^{2}\widetilde{\psi}^{2}_{s}(u)\mu(du)ds}\Bigr]
\leq
E\Bigl[\sqrt{ess\sup_{\tau\in\mathcal{T}_{0,T}}\widetilde{Y}_{\tau}^{2}}\sqrt{\int_{0}^{T}e^{2\beta
s}\int_{\mathcal{U}}\widetilde{\psi}^{2}_{s}(u)\mu(du)ds}\Bigr]$$
Thus
\begin{equation}
  E\Bigl[\sqrt{\int_{0}^{T}\int_{\mathcal{U}}e^{2\beta
s}\widetilde{Y}_{s-}^{2}\widetilde{\psi}^{2}_{s}(u)\mu(du)ds}\Bigr]
\leq
\||\widetilde{Y}\||_{\mathds{S}^{2}}\|\widetilde{\psi}\|_{\mathbb{L}^{2,2\beta}_{\pi}}<\infty
   \label{moneq17}
\end{equation}
Then  $E\Bigl[\int_{0}^{T}\int_{\mathcal{U}}e^{\beta
s}\widetilde{Y}_{s-}\widetilde{\psi}_{s}(u)\widetilde{\pi}(du,ds)\Bigr]=0$.
Finally the same result holds for the martingale
$\int_{0}^{t}e^{\beta s}\widetilde{Y}_{s-}d\widetilde{M}_{s}$,
since:
\begin{equation}
E\Bigl[\sqrt{\int_{0}^{T} e^{2\beta
s}\widetilde{Y}_{s-}^{2}d[\widetilde{M}]_{s}}\Bigr] \leq
\||\widetilde{Y}\||_{\mathds{S}^{2}}
\|\widetilde{M}\|_{\mathbb{M}_{2\beta}^{2}} <\infty \label{moneq18}
\end{equation}
By taking expectations on both sides of $~\eqref{moneq11}$ with
$t=0$, we obtain: $$\widetilde{Y}^{2}_{0}+E[\int_{0}^{T} e^{\beta
s}\widetilde{Z}^{2}_{s}ds]+E[\int_{0}^{T} e^{\beta
s}d<\widetilde{M}>_{s}]+E[\int_{0}^{T} e^{\beta
s}\int_{\mathcal{U}}|\widetilde{\psi}_{s}(u)|^{2}\mu(du)ds]\leq
\epsilon^{2}\|\widetilde{f}(s)\|^{2}_{\beta}.$$ Hence, with the fact
that  $E[\int_{0}^{T} e^{\beta
s}d<\widetilde{M}>_{s}]=E[\int_{0}^{T} e^{\beta
s}d[\widetilde{M}]_{s}]$, we have

\begin{equation}
\|\widetilde{Z}\|^{2}_{\beta}+\|\widetilde{M}\|_{\mathbb{M}_{\beta}^{2}}^{2}
    +\|\widetilde{\psi}\|^{2}_{\mathbb{L}^{2,\beta}_{\pi}}\leq \epsilon^{2}\|\widetilde{f}(s)\|^{2}_{\beta}
\label{moneq19}
\end{equation}
This therefore shows the first inequality of the lemma. From
$~\eqref{moneq11}$ we also get, for all $\tau\in\mathcal{T}_{0,T}$
\begin{multline}
   e^{\beta\tau}\widetilde{Y}_{\tau}^{2}\leq
    \epsilon^{2}\int_{0}^{T}e^{\beta s}\widetilde{f}^{2}(s)ds-2\int_{\tau}^{T}e^{\beta s}\widetilde{Y}_{s-}\widetilde{Z}_{s}dW_{s}
     -2\int_{\tau}^{T}\int_{\mathcal{U}}e^{\beta s}\widetilde{Y}_{s-}\widetilde{\psi}_{s}(u)\widetilde{\pi}(du,ds)\\ -
     2\int_{\tau}^{T}e^{\beta s}\widetilde{Y}_{s-}d\widetilde{M}_{s} \,\,\ a.s.
 \label{moneq20}
\end{multline}
By taking first the essential supremum over
$\tau\in\mathcal{T}_{0,T}$, and then the expectation on both sides
of the  inequality $~\eqref{moneq20}$, we obtain:
\begin{multline}
E\Bigl[ess\sup_{\tau\in\mathcal{T}_{0,T}}e^{\beta
\tau}\widetilde{Y}_{\tau}^{2}\Bigr] \leq
\epsilon^{2}\|\widetilde{f}(s)\|^{2}_{\beta}+2E\Bigl[ess\sup_{\tau\in\mathcal{T}_{0,T}}|\int_{0}^{\tau}e^{\beta s}\widetilde{Y}_{s-}\widetilde{Z}_{s}dW_{s}|\Bigr]\\
+2E\Bigl[ess\sup_{\tau\in\mathcal{T}_{0,T}}|\int_{0}^{\tau}\int_{\mathcal{U}}e^{\beta s}\widetilde{Y}_{s-}\widetilde{\psi}_{s}(u)\widetilde{\pi}(du,ds)|\Bigr]\\
+2E\Bigl[ess\sup_{\tau\in\mathcal{T}_{0,T}}|\int_{0}^{\tau}e^{\beta
s}\widetilde{Y}_{s-}d\widetilde{M}_{s}|\Bigr].
 \label{moneq21}
\end{multline}
By using the continuity of a.e. trajectory of the process
$(\int_{0}^{t}e^{\beta
s}\widetilde{Y}_{s-}\widetilde{Z}_{s}dW_{s})_{t\in[0,T]}$
(\cite{Ouknine2015}, Prop.A.3 ) and Burkholder-Davis-Gundy
inequalities (\cite{Protter2000} Theorem 48, page 193. Applied
with $p=1$), we get
\begin{multline}
E\Bigl[ess\sup_{\tau\in\mathcal{T}_{0,T}}|\int_{0}^{\tau} e^{\beta
s}\widetilde{Y}_{s-}\widetilde{Z}_{s}dW_{s}|\Bigr]=E\Bigl[\sup_{t\in[0,T]}|\int_{0}^{t}
e^{\beta s}\widetilde{Y}_{s-}\widetilde{Z}_{s}dW_{s}|\Bigr]\\
\leq cE\Bigl[\sqrt{\int_{0}^{T}e^{2\beta s}\widetilde{Y}_{s-}^{2}\widetilde{Z}^{2}_{s}ds} \Bigr]
\label{moneq22}
\end{multline}
where $c$ is a positive "universal" constant (which does not depend
on the other parameters). The same reasoning as that used to obtain
equation $~\eqref{moneq14}$ leads to
\begin{equation}
\sqrt{\int_{0}^{T}e^{2\beta
s}\widetilde{Y}_{s-}^{2}\widetilde{Z}^{2}_{s}ds}\leq
\sqrt{ess\sup_{\tau\in\mathcal{T}_{0,T}}e^{\beta
\tau}\widetilde{Y}_{\tau}^{2}\int_{0}^{T}e^{\beta
s}\widetilde{Z}^{2}_{s}ds} \,\,\ p.s. \label{moneq23}
\end{equation}
From the  inequalities  $~\eqref{moneq22}$, $~\eqref{moneq23}$ and
$ab\leq\frac{1}{4}a^{2}+b^{2}$, we have
\begin{equation}
2E\Bigl[ess\sup_{\tau\in\mathcal{T}_{0,T}}|\int_{0}^{\tau} e^{\beta
s}\widetilde{Y}_{s-}\widetilde{Z}_{s}ds|\Bigr]\leq\frac{1}{4}E\Bigl[ess\sup_{\tau\in\mathcal{T}_{0,T}}e^{\beta
\tau}\widetilde{Y}_{\tau}^{2}\Bigr]+4c^{2}E\Bigl[\int_{0}^{T}e^{\beta
s}\widetilde{Z}^{2}_{s}ds\Bigr]. \label{moneq24}
\end{equation}
By the same arguments, we have
\begin{multline}
2E\Bigl[ess\sup_{\tau\in\mathcal{T}_{0,T}}|\int_{0}^{\tau}e^{\beta
s}\widetilde{Y}_{s-}\widetilde{\psi}_{s}(u)\widetilde{\pi}(du,ds)|\Bigr]\leq
2cE\Bigl[\sqrt{|\int_{0}^{T}\int_{\mathcal{U}}e^{2\beta s}\widetilde{Y}^{2}_{s-}\widetilde{\psi}^{2}_{s}(u)\mu(du)ds|}\Bigr]\\
\leq \frac{1}{4}E\Bigl[ess\sup_{\tau\in\mathcal{T}_{0,T}}e^{\beta
\tau}\widetilde{Y}_{\tau}^{2}\Bigr]+4c^{2}E\Bigl[\int_{0}^{T}\int_{\mathcal{U}}e^{\beta
s}\widetilde{\psi}^{2}_{s}(u)\mu(du)ds\Bigr] \label{moneq25}
\end{multline}
And
\begin{multline}
2E\Bigl[ess\sup_{\tau\in\mathcal{T}_{0,T}}|\int_{0}^{\tau}e^{\beta
s}\widetilde{Y}_{s-}d\widetilde{M}_{s}|\Bigr]\leq
2cE\Bigl[\sqrt{|\int_{0}^{T}e^{2\beta s}\widetilde{Y}^{2}_{s-}d[\widetilde{M}]_{s}|}\Bigr]\\
\leq \frac{1}{4}E\Bigl[ess\sup_{\tau\in\mathcal{T}_{0,T}}e^{\beta
\tau}\widetilde{Y}_{\tau}^{2}\Bigr]+4c^{2}E\Bigl[\int_{0}^{T}e^{\beta
s}d[\widetilde{M}]_{s}\Bigr] \label{moneq26}
\end{multline}
where $c$ is a positive  constant which does not depend on the other
parameters. From $~\eqref{moneq22}$, $~\eqref{moneq24}$,
$~\eqref{moneq25}$ and  $~\eqref{moneq26}$, we get
$$\frac{1}{4}\||\widetilde{Y}|\|_{\beta}^{2}\leq\epsilon^{2}\|\widetilde{f}(s)\|^{2}_{\beta}
+4c^{2}\|\widetilde{Z}\|^{2}_{\beta}
+4c^{2}\|\widetilde{M}\|_{\mathbb{M}_{\beta}^{2}}^{2}
+4c^{2}\|\widetilde{\psi}\|^{2}_{\mathbb{L}^{2,\beta}_{\pi}}$$ This
inequality, combined with  $~\eqref{moneq19}$, gives
$$\||\widetilde{Y}|\|_{\beta}^{2}\leq 4\epsilon^{2}(1+4c^{2})\|\widetilde{f}(s)\|^{2}_{\beta}$$
\end{proof}
In the following lemma, we prove existence and uniqueness of the
solution to the RBSDE from Definition$~\ref{definition1}$ in the
case where the driver $f$ depends only on $s$ and $\omega$, i.e.
$f(\omega,s,y,z,\psi):=f(\omega,s)$.
\begin{lemma}
\label{lemme3} Suppose that $f$ does not depend on y, z, $\psi$ that
is $f(\omega,s,y,z,\psi):=f(\omega,s)$, where $f$ is a process in
$\mathds{H}^{2}$. Let $\xi$ be an obstacle. Then, the RBSDE from
Definition$~\ref{definition1}$   admits a unique solution
$(Y,Z,\psi,M,A,C)\in\mathcal{E}^{2}(0,T)$, and for each
$S\in\mathcal{T}_{0,T}$, we have
\begin{equation}
Y_{S}=ess\sup_{\tau\in\mathcal{T}_{S,T}}E\Bigl[\xi_{\tau}+\int_{S}^{\tau}
f(t)dt|\mathcal{F}_{S} \Bigr] \,\,\,\ a.s. \label{moneq27}
\end{equation}
\end{lemma}

\begin{proof}
For all  $S\in\mathcal{T}_{0,T}$, we define $\overline{Y}(S)$  by:
\begin{equation}
\overline{Y}(S)=ess\sup_{\tau\in\mathcal{T}_{S,T}}E\Bigl[\xi_{\tau}+\int_{S}^{\tau}
f(t)dt|\mathcal{F}_{S} \Bigr] \,\,\,\  , \,\,\,\
\overline{Y}(T)=\xi_{T} \label{moneq28}
\end{equation}
And $\overline{\overline{Y}}(S)$  by:
\begin{equation}
\overline{\overline{Y}}(S)=\overline{Y}(S)+\int_{0}^{S}f(t)dt=ess\sup_{\tau\in\mathcal{T}_{S,T}}E\Bigl[\xi_{\tau}+\int_{0}^{\tau}
f(t)dt|\mathcal{F}_{S} \Bigr]
 \label{moneq29}
\end{equation}
We note that the process $(\xi_{t}+\int_{0}^{t} f(s)ds)_{t\in[0,T]}$
is progressive. Therefore, the family
$(\overline{\overline{Y}}(S))_{S\in\mathcal{T}_{0,T}}$ is a
supermartingale family (see \cite{KobylanskiQuenez2012} Remark
1.2 with Prop.1.5), and with remark $(b)$ in
(\cite{DellacherieMeyer1980}, page 435), gives the existence of a
strong optional supermartingale (which we denote again by
$\overline{\overline{Y}}$) such that
$\overline{\overline{Y}}_{S}=\overline{\overline{Y}}(S)$ a.s. for
all $S\in\mathcal{T}_{0,T}$. Thus, we have
$\overline{Y}(S)=\overline{\overline{Y}}(S)-\int_{0}^{S}f(t)dt=\overline{\overline{Y}}_{S}-\int_{0}^{S}f(t)dt$
a.s. for all  $S\in\mathcal{T}_{0,T}$ (see
\cite{DellacherieMeyer1980}). On the other hand, we know that
almost all trajectories of the strong optional supermartingale
$\overline{\overline{Y}}$  are l\`{a}dl\`{a}g. Thus, we get that the l\`{a}dl\`{a}g
optional process
$(\overline{Y}_{t})_{t\in[0,T]}=(\overline{\overline{Y}}_{t}-\int_{0}^{t}f(s)ds)_{t\in[0,T]}$
aggregates the family $(\overline{Y}(S))_{S\mathcal{T}_{0,T}}$.

To prove the lemma $~\ref{lemme3}$, it must be shown, as a first
step, that $\overline{Y}\in\mathds{S}^{2}$ by giving  an estimate of
$\||\overline{Y}|\|_{\mathds{S}^{2}}^{2}$ in terms of
$\||\xi|\|_{\mathds{S}^{2}}^{2}$ and $\|f\|_{\mathds{H}^{2}}^{2}$.
In the second step, we exhibit processes $Z$, $\psi$, $M$, $A$ and
$C$ such that $(\overline{Y},Z,\psi,M,A,C)$ is a solution to the
RBSDE with parameters $(f,\xi)$. In the third step, we prove that
$A\times C\in \mathds{S}^{2}\times\mathds{S}^{2}$  and we give an
estimate of $\||A|\|_{\mathds{S}^{2}}^{2}$ and
$\||C|\|_{\mathds{S}^{2}}^{2}$. In the fourth step, we show that
$Z\in\mathds{H}^{2}$, $\psi\in\mathbb{L}_{\pi}^{2}$ and
$M\in\mathbb{M}^{2}$, and finally we show the uniqueness of the
solution.\\
\textbf{Step 1.} By using the definition of $\overline{Y}$
$~\eqref{moneq28}$, Jensen's inequality and the triangular
inequality, we get
$$|\overline{Y}_{S}|\leq ess\sup_{\tau\in\mathcal{T}_{S,T}}E\Bigl[|\xi_{\tau}|+|\int_{S}^{\tau}
f(t)dt| |\mathcal{F}_{S} \Bigr]\leq
E\Bigl[ess\sup_{\tau\in\mathcal{T}_{S,T}}|\xi_{\tau}|+\int_{0}^{T}
|f(t)|dt |\mathcal{F}_{S} \Bigr]$$ Thus, we obtain
\begin{equation}
|\overline{Y}_{S}|\leq E\Bigl[X|\mathcal{F}_{S}\Bigr]
 \label{moneq30}
\end{equation}
With
\begin{equation}
X=\int_{0}^{T}|f(t)|dt+ess\sup_{\tau\in\mathcal{T}_{0,T}}|\xi_{\tau}|
 \label{moneq31}
\end{equation}
Applying Cauchy-Schwarz inequality gives
\begin{equation}
E[X^2]\leq
cT\|f\|_{\mathds{H}^{2}}^{2}+c\||\xi|\|_{\mathds{S}^{2}}^{2}<\infty.
\label{moneq32}
\end{equation}
where $c$ is a positive constant. Now, inequality $~\eqref{moneq30}$
leads to $|\overline{Y}_{S}|^{2}\leq |E[X|\mathcal{F}_{S}]|^{2}$. By
taking the essential supremum over $S\in\mathcal{T}_{0,T}$ we get
$ess\sup_{S\in\mathcal{T}_{0,T}}|\overline{Y}_{S}|^{2}\leq
ess\sup_{S\in\mathcal{T}_{0,T}}|E[X|\mathcal{F}_{S}]|^{2}$. By using
Proposition A.3 in \cite{Ouknine2015}, we get
$ess\sup_{S\in\mathcal{T}_{0,T}}|\overline{Y}_{S}|^{2}\leq
\sup_{t\in[0,T]}|E[X|\mathcal{F}_{t}]|^{2}$. By using this
inequality and Doob's martingale inequalities, we obtain
\begin{equation}
E\Bigl[ess\sup_{S\in\mathcal{T}_{0,T}}|\overline{Y}_{S}|^{2}\Bigr]\leq
E\Bigl[\sup_{t\in[0,T]}|E[X|\mathcal{F}_{t}]|^{2}\Bigr] \leq
cE[X^{2}] \label{moneq33}
\end{equation}
where $c$ is a positive constant  that changes from line to line.
Finally, combining inequalities $~\eqref{moneq32}$ and
$~\eqref{moneq33}$ gives
\begin{equation}
E\Bigl[ess\sup_{S\in\mathcal{T}_{0,T}}|\overline{Y}_{S}|^{2}\Bigr]\leq
cT\|f\|_{\mathds{H}^{2}}^{2}+c\||\xi|\|_{\mathds{S}^{2}}^{2}<\infty.
\label{moneq34}
\end{equation}
Then $\overline{Y}_{S}\in\mathds{S}^{2}$. \\
\textbf{Step 2.} Due to the previous step and to the assumption
$f\in\mathds{H}^{2}$, the strong optional supermartingale
$\overline{\overline{Y}}$ is of class $(D)$. Applying Mertens
decomposition (\cite{Ouknine2015}, Theorem A.1) and a result from
optimal stopping theory (see more in \cite{ElKaroui1981}, Prop.
2.34. page 131 or \cite{KobylanskiQuenez2012}), gives the
following
$$\overline{\overline{Y}}_{\tau}=N_{\tau}-A_{\tau}-C_{\tau-} \,\,\,\ \forall \tau\in\mathcal{T}_{0,T}$$
\begin{equation}
\overline{Y}_{\tau}=-\int_{0}^{\tau}f(t)dt+N_{\tau}-A_{\tau}-C_{\tau-}
\,\,\,\ a.s. \,\,\,\ \forall \tau\in\mathcal{T}_{0,T}
 \label{moneq35}
\end{equation}
where $N$ is a (c\`{a}dl\`{a}g) uniformly integrable martingale such that
$N_0=0$, $A$ is a nondecreasing right-continuous predictable process
such that $A_0=0$, $E(A_T)<\infty$ and satisfying $~\eqref{moneq6}$,
and $C$ is a nondecreasing right-continuous adapted purely
discontinuous process such that $C_{0-}=0$, $E(C_T)<\infty$ and
satisfying $~\eqref{moneq7}$. By the martingale representation
theorem (Lemma$~\ref{lemme1}$), there exists a unique predictable
process $Z$, a unique  process $\psi$ and a unique (c\`{a}dl\`{a}g) local
martingales orthogonal $M$  such that
$$N_{t}=\int_{0}^{t}Z_{s}dW_{s}+\int_{0}^{t}\int_{\mathcal{U}}\psi_{s}(u)\widetilde{\pi}(du,ds)+M_{t}.$$
Moreover, we have $\overline{Y}_{T}=\xi_{T}$ a.s. by definition of
$\overline{Y}$. Combining this with equation $~\eqref{moneq35}$.
gives equation $~\eqref{moneq4}$. Also by definition of
$\overline{Y}$, we have $\overline{Y}_{S}\geq\xi_{S}$ a.s. for all
$S\in\mathcal{T}_{0,T}$, which, along with Proposition A.4 in
\cite{Ouknine2015} (or Theorem 3.2. in
\cite{AshkanNikeghbali2006}), shows that $\overline{Y}$ satisfies
inequality $~\eqref{moneq5}$. Finally, to conclude that the process
$(\overline{Y},Z,\psi,M,A,C)$ is a solution to the RBSDE with
parameters $(f,\xi)$, it remains to show that $Z\times\psi\times
M\times A\times
C\in\mathds{H}^{2}\times\mathbb{L}_{\pi}^{2}\times\mathbb{M}^{2}\times\mathds{S}^{2}\times\mathds{S}^{2}$.\\
\textbf{Step 3.} Let us show that $A\times C\in\mathds{S}^{2}\times\mathds{S}^{2}$.\\
Let us define the process $\overline{\overline{A}}_{t}=A_{t}+C_{t-}$
where the processes $A$ and  $C$ are given by $~\eqref{moneq35}$. By
arguments similar to those used in the proof of inequality
$~\eqref{moneq30}$, we see that $|\overline{\overline{Y}}_{S}|\leq
E[X|\mathcal{F}_{S}]$ with
$$X=\int_{0}^{T}|f(t)|dt+ess\sup_{\tau\in\mathcal{T}_{S,T}}|\xi_{\tau}|.$$
Then, the Corollary A.1 in \cite{Ouknine2015} ensures the
existence of a constant $c>0$ such that
$E[(\overline{\overline{A}}_{T})^{2}]\leq cE[X^{2}]$. By combining
this inequality with inequality $~\eqref{moneq32}$, we obtain
\begin{equation}
E[(\overline{\overline{A}}_{T})^{2}]\leq
cT\|f\|_{\mathds{H}^{2}}^{2}+c\||\xi|\|_{\mathds{S}^{2}}^{2}
 \label{moneq36}
\end{equation}
where we have again allowed the positive constant $c$ to vary from
line to line. We conclude that  $\overline{\overline{A}}\in L^{2}$.
And with the nondecreasingness of $\overline{\overline{A}}$, then
$(\overline{\overline{A}}_{\tau})^{2}\leq
(\overline{\overline{A}}_{T})^{2}$ for all
$\tau\in\mathcal{T}_{0,T}$ thus
$$E\Bigl[ess\sup_{\tau\in\mathcal{T}_{0,T}}(\overline{\overline{A}}_{\tau})^{2}\Bigr]\leq
E\Bigl[(\overline{\overline{A}}_{T})^{2}\Bigr]$$ i.e.
$\overline{\overline{A}}\in \mathds{S}^{2}$ then $A\in
\mathds{S}^{2}$ and $C\in \mathds{S}^{2}$. \\
\textbf{Step 4.} Let us now prove that $Z\times\psi\times
M\in\mathds{H}^{2}\times\mathbb{L}_{\pi}^{2} \times\mathbb{M}^{2}$.
We have from step 3
$$\int_{0}^{T}Z_{s}dW_{s}+\int_{0}^{T}\int_{\mathcal{U}}\psi_{s}(u)\widetilde{\pi}(du,ds)+M_{T}=
\overline{Y}_{T}+\int_{0}^{T}f(t)dt+\overline{\overline{A}}_{T}-\overline{Y}_{0}$$
where $\overline{\overline{A}}$ is the process from Step 3. Since
$\overline{\overline{A}}_{T}\in L^{2}$, $\overline{Y}_{T}\in L^{2}$,
 $\overline{Y}_{0}\in L^{2}$ and $f\in\mathds{H}^{2}$. Hence,  $\int_{0}^{T}Z_{s}dW_{s}\in L^{2}$,
 $\int_{0}^{T}\int_{\mathcal{U}}\psi_{s}(u)\widetilde{\pi}(du,ds)\in L^{2}$ and $M_{T}\in L^{2}$ and consequently
$Z\times\psi\times M\in\mathds{H}^{2}\times\mathbb{L}_{\pi}^{2}
\times\mathbb{M}^{2}$.

For the uniqueness of the solution, suppose that $(Y,Z,\psi,M,A,C)$
is a solution of the RBSDE with driver $f$ and obstacle $\xi$. Then,
by the previous inequality$~\ref{moneq8}$ in the
Lemma$~\ref{lemme2}$ (applied with $f^1=f^2=f$) we obtain
$Y=\overline{Y}$ in $\mathds{S}^{2}$, where $\overline{Y}$ is given
by $~\eqref{moneq28}$. The uniqueness of $A$, $C$, $Z$, $\psi$ and
$M$ follows from the uniqueness of Mertens decomposition of strong
optional supermartingales and from the uniqueness of the martingale
representation (Lemma$~\ref{lemme1}$).
\end{proof}

\begin{remark}
\label{remarque1}
\begin{enumerate}
    \item We note that the uniqueness of $Z$, $\psi$  and $M$ can be obtained also by applying $~\eqref{unicite}$ in the previous Lemma$~\ref{lemme2}$.
    \item Let $\beta>0$. For $\phi\in\mathds{S}^{2}$, we have
    $E[\int_{0}^{T}e^{\beta t}|\phi_{t}|^{2}dt]\leq TE[ess\sup_{\tau\in\mathcal{T}_{0,T}}e^{\beta \tau}|\phi_{\tau}|^{2}]$.
    Indeed, by applying Fubini's theorem, we get
\begin{multline}
E[\int_{0}^{T}e^{\beta t}|\phi_{t}|^{2}dt]= \int_{0}^{T}E[e^{\beta
t}|\phi_{t}|^{2}]dt\leq
\int_{0}^{T}E[ess\sup_{\tau\in\mathcal{T}_{0,T}}e^{\beta \tau}|\phi_{\tau}|^{2}]ds=\\
TE[ess\sup_{\tau\in\mathcal{T}_{0,T}}e^{\beta
\tau}|\phi_{\tau}|^{2}]
\end{multline}
\end{enumerate}
\end{remark}
In the following theorem, we prove existence and uniqueness of the
solution to the RBSDE from Definition$~\ref{definition1}$ in the
case of a general Lipschitz driver $f$ by using a fixed-point
theorem and by using (2) in the Remark$~\ref{remarque1}$ .
\begin{theorem}
Let $\xi$ be a left-limited and r.u.s.c. process in $\mathds{S}^{2}$
and let $f$ be a Lipschitz driver. The RBSDE with parameters
$(f,\xi)$ from Definition$~\ref{definition1}$  admits a unique
solution $(Y,Z,\psi,M,A,C)\in\mathcal{E}^{2}(0,T)$.\\

\end{theorem}
\begin{proof}
We note by $\mathcal{E}^{\beta}_{f}$ the space
$\mathds{S}^{2}\times\mathds{H}^{2}\times\mathbb{L}^{2}_{\pi}(0,T)$
which we equip with the norm $\|.\|_{\mathcal{E}^{\beta}_{f}}$
defined by
$$\|(Y,Z,\psi)\|_{\mathcal{E}^{\beta}_{f}}^{2}=\||Y|\|_{\beta}^{2}+\|Z\|^{2}_{\beta}+\|\psi\|^{2}_{\mathbb{L}^{2,\beta}_{\pi}}$$ for all
$(Y,Z,\psi)\in\mathds{S}^{2}\times\mathds{H}^{2}\times\mathbb{L}^{2}_{\pi}(0,T)$.
After, we define an application
$\Phi:\mathcal{E}^{\beta}_{f}\rightarrow \mathcal{E}^{\beta}_{f}$ as
follows: for a given $(y,z,\varphi)\in\mathcal{E}^{\beta}_{f}$, we
let $(Y,Z,\psi)=\Phi(y,z,\varphi)$ where $(Y,Z,\psi)$ the first
three components of the solution to the RBSDE associated with driver
$f:=f(t,y_t,z_t,\varphi_t)$ and with obstacle $\xi_t$. Let $(A,C)$
be the associated Mertens process, constructed as in
lemma$~\ref{lemme3}$. The mapping $\Phi$ is well-defined by
Lemma$~\ref{lemme3}$.

Let $(y,z,\varphi)$ and $(y',z',\varphi')$ be two elements of
$\mathcal{E}^{\beta}_{f}$. We set $(Y,Z,\psi)=\Phi(y,z,\varphi)$
 and  $(Y',Z',\psi')=\Phi(y',z',\varphi')$. We also set $\widetilde{Y}=Y-Y'$, $\widetilde{Z}=Z-Z'$,
$\widetilde{\psi}=\psi-\psi'$, $\widetilde{y}=y-y'$,
$\widetilde{z}=z-z'$ and $\widetilde{\varphi}=\varphi-\varphi'$.

By the same argument that in the proof of Theorem 3.4 in
\cite{Ouknine2015}, in the Brownian filtration case. Let us prove
that for a suitable choice of the parameter $\beta>0$, the mapping
$\Phi$ is a contraction from the Banach space
$\mathcal{E}^{\beta}_{f}$ into itself. Indeed, By applying
Lemma$~\ref{lemme2}$, we have, for all $\epsilon>0$ and for all
$\beta\geq\frac{1}{\epsilon^{2}}$:
\begin{eqnarray*}
\||\widetilde{Y}|\|_{\beta}^{2}+\|\widetilde{Z}\|^{2}_{\beta}+\|\widetilde{\psi}\|^{2}_{\mathbb{L}^{2,\beta}_{\pi}}&\leq
&\||\widetilde{Y}|\|_{\beta}^{2}+\|\widetilde{Z}\|^{2}_{\beta}+\|\widetilde{M}\|^{2}_{\mathbb{M}^{2}_{\beta}}
    +\|\widetilde{\psi}\|^{2}_{\mathbb{L}^{2,\beta}_{\pi}}\\
&\leq&
\epsilon^{2}(5+16c^{2})\|f(t,y,z,\varphi)-f(t,y',z',\varphi')\|_{\beta}^{2}
\end{eqnarray*}
By using the Lipschitz property of $f$ and the fact that
$(a+b)^2\leq2a^2+2b^2$, for all $(a,b)\in\R^{2}$, we obtain
 $$\|f(t,y,z,\varphi)-f(t,y',z',\varphi')\|_{\beta}^{2}\leq
C_{K}(\|\widetilde{y}\|_{\beta}^{2}+\|\widetilde{z}\|_{\beta}^{2}
+\|\widetilde{\varphi}\|^{2}_{\mathbb{L}^{2,\beta}_{\pi}})$$ where
$C_{K}$ is a positive constant depending on the Lipschitz constant
$K$ only. Thus, for all $\epsilon>0$ and for all
$\beta\geq\frac{1}{\epsilon^{2}}$ we have:
$$\||\widetilde{Y}|\|_{\beta}^{2}+\|\widetilde{Z}\|^{2}_{\beta}+\|\widetilde{\psi}\|^{2}_{\mathbb{L}^{2,\beta}_{\pi}}\leq
    \epsilon^{2}C_{K}(5+16c^{2})\Bigl(\|\widetilde{y}\|_{\beta}^{2}+\|\widetilde{z}\|_{\beta}^{2}
+\|\widetilde{\varphi}\|^{2}_{\mathbb{L}^{2,\beta}_{\pi}}\Bigr)$$
The previous inequality, combined with (2) in
Remark$~\ref{remarque1}$, gives
$$\||\widetilde{Y}|\|_{\beta}^{2}+\|\widetilde{Z}\|^{2}_{\beta}+\|\widetilde{\psi}\|^{2}_{\mathbb{L}^{2,\beta}_{\pi}}\leq
    \epsilon^{2}C_{K}(5+16c^{2})(T+1)\Bigl(\||\widetilde{y}|\|_{\beta}^{2}+\|\widetilde{z}\|_{\beta}^{2}
+\|\widetilde{\varphi}\|^{2}_{\mathbb{L}^{2,\beta}_{\pi}}\Bigr)$$
Thus, for $\epsilon>0$ such that
$\epsilon^{2}C_{K}(5+16c^{2})(T+1)<1$ and $\beta>0$ such that
$\beta\geq\frac{1}{\epsilon^{2}}$, the mapping $\Phi$ is a
contraction. By the Banach fixed-point theorem, we get that $\Phi$
has a unique fixed point in $\mathcal{E}^{\beta}_{f}$. We thus have
the existence and uniqueness of the solution to the RBSDE.
\end{proof}
\section{Application on dynamic risk measure and optimal stopping.}
\subsection{On dynamic risk measure}
In this subsection, we give  an application of reflected BSDEs in dynamic risk measure. Indeed, define the following functional: for each stopping time $\tau\in\mathcal{T}_{0,T}$ and $\xi\in\mathds{S}^{2}$. Set
\begin{equation}
v(S)=-ess\sup_{\tau\in\mathcal{T}_{S,T}} \mathcal{E}^{f}_{S,\tau}(\xi_{\tau})
 \label{moneq38}
\end{equation}
where $S\in\mathcal{T}_{0,T}$, $v$ is the dynamic risk measure, $\xi_{T'}$ $(T'\in[0,T])$ is the gain of the position at time $T'$ and $-\mathcal{E}^{f}_{t,T'}(\xi_{T'})$ is the $f$-conditional expectation of $\xi_{\tau}$ modelling the risk at time $t$ where $t\in[0,T]$. We can show that the minimal risk measure $v$ defined by $~\eqref{moneq38}$ coincides with $-Y$, where
$Y$ is (the first component of) the solution to the reflected BSDE associated with driver $f$ and obstacle
$\xi$. For this purpose, we can extend the results in Proposition A.5 and Theorem 4.2 in \cite{Ouknine2015} to our setting (see \cite{AaziziOuknine2016} for more details).
\subsection{On optimal stopping}
We note also that we can show the existence of an $\varepsilon$-optimal stopping time, and that of the
existence of an optimal stopping time under suitable assumptions on the barrier  $\xi$ i.e. without right continuity of $\xi$, by extending the results of the second part of \cite{Ouknine2015} to our setting.

Let $(Y,Z,\psi,M,A,C)$ be the solution of the reflected BSDE with parameters $(f,\xi)$ as in definition$~\ref{definition1}$, we have
\begin{equation}
Y_{S}=ess\sup_{\tau\in\mathcal{T}_{S,T}} \mathcal{E}^{f}_{S,\tau}(\xi_{\tau})
 \label{moneq39}
\end{equation}
For each $S\in\mathcal{T}_{0,T}$ and $\varepsilon >0$, the stopping time $\tau_{S}^{\varepsilon}=inf \{t\geq S, Y_{t}\leq\xi_{t}+\varepsilon\}$ is a $(C\varepsilon)$-optimal for $~\ref{moneq39}$ where $C$ is a constant which  depends only on $T$  and the Lipschitz constant $K$ of $f$:
$$Y_{S}\leq\mathcal{E}^{f}_{S,\tau_{S}^{\varepsilon}}(\xi_{\tau_{S}^{\varepsilon}})+C\varepsilon, \,\,\,\ a.s.$$
Under our assumption on $\xi$ and $f$, we can prove that for each $S\in\mathcal{T}_{0,T}$ and $\widehat{\tau}\in\mathcal{T}_{S,T}$, the stopping time $\widehat{\tau}$ is $S$-optimal. i.e. $$Y_{\widehat{\tau}}=\mathcal{E}^{f}_{S,\widehat{\tau}}(\xi_{\widehat{\tau}}), \,\,\,\ a.s.$$
 (see Theorem 4.2 and Proposition 4.3 in \cite{Ouknine2015}).

 Finally, under an additional assumption of left-upper semicontinuity (l.u.s.c) of $\xi$ in $\mathds{S}^{2}$, the first time when the value process $Y$ hits $\xi$ is optimal: if  $\tau_{S}^{*}=inf \{u\geq S, Y_{u}=\xi_{u}\}$, $\xi$ is r.u.s.c and l.u.s.c in $\mathds{S}^{2}$ and  $(Y,Z,\psi,M,A,C)$ is the solution of the reflected BSDE of definition$~\ref{definition1}$, the stopping time  $\tau_{S}^{*}$ is optimal that is
 $$Y_{S}=\mathcal{E}^{f}_{S,\tau_{S}^{*}}(\xi_{\tau_{S}^{*}}), \,\,\,\ a.s.$$
 (see Proposition 4.2 in \cite{Ouknine2015}).

\section*{Acknowledgements}
The authors would like to thank the referee for the careful reading of the paper and highly appreciate the comments and suggestions, which significantly contributed to improving the quality of the paper.

\bibliographystyle{alea3}
\bibliography{ALEArbsde}

\end{document}